\documentclass[12pt]{article}

\usepackage[dvips]{graphicx}
\usepackage{t1enc,amssymb,amsbsy,amsthm,amsmath,amsfonts}

\newcommand{\rr}{\stackrel {d}{=}}

\renewcommand{\Re}{{\rm I\kern-0.16em R}}

\def\@begintheorem#1#2{\trivlist \item[\hskip \labelsep{\bf #1\ #2}]}
\def\@opargbegintheorem#1#2#3{\trivlist
      \item[\hskip \labelsep{\bf #1\ #2\ (#3)}]}

\newtheorem{proposition}{Proposition}[section] 

\newtheorem{lemma}[proposition]{Lemma}
\newtheorem{corollary}[proposition]{Corollary}
\newtheorem{example}[proposition]{Example}
\newtheorem{remark}[proposition]{Remark}

\def\R{{\bf R}}
\def\R{{\bf R}}

\def\E{{\bf E}}

\def\cF{{\cal F}}

\def\cB{{\cal B}}

\def\cG{{\cal G}}

\def\cK{{\cal K}}

\def\al{\alpha}

\def\ga{\gamma}

\numberwithin{equation}{section}

\begin{document}

\author{Paavo Salminen\\{\small Åbo Akademi,}
\\{\small Mathematical Department,}
\\{\small FIN-20500 Åbo, Finland,} \\
{\small email: phsalmin@abo.fi}
\and
Marc Yor\\
{\small Universit\'e Pierre et Marie Curie,}\\
{\small Laboratoire de Probabilit\'es }\\
{\small et Mod\`eles al\'eatoires ,}\\
{\small 4, Place Jussieu, Case 188}\\
{\small F-75252 Paris Cedex 05, France}\\
}
\vskip5cm


\title{Tanaka formula for symmetric L\'evy processes}

\maketitle

\begin{abstract}
Starting from the potential  theoretic definition of the local times
of a Markov process -- when these exist -- we obtain a Tanaka formula
for the local times of symmetric L\'evy processes. The most interesting
case is that of the symmetric $\al$-stable L\'evy process (for $\al\in(1,2]$)
which is studied in detail. In particular, we determine which powers
of such a process are semimartingales. These results complete, in a
sense, the works by K. Yamada \cite{yamada02} and Fitzsimmons and
Getoor \cite{fitzsimmonsgetoor92a}. 
\end{abstract}

\section{Introduction and main results}
\label{sec1}
It is
well known that there are different constructions
and definitions of local times corresponding to different classes of
stochastic processes. For a large panorama of such definitions, see
Geman and Horowitz \cite{gemanhorowitz80}.

The most common definition of the local times $L=\{L^x_t\,:\,x\in\R ,\, t\geq
0\}$ of a given process $\{X_t\,:\, t\geq 0\}$ is as the
Radon--Nikodym derivative of the occupation measure of $X$ with
respect to the Lebesgue measure in $\R;$ precisely $L$
satisfies
\begin{equation}     
\label{eq1}
\int_0^t f(X_s)\,ds =
\int_{-\infty}^\infty f(x) L^x_t\,dx
\end{equation}
for every Borel function $f:\R\mapsto\R_+ .$

There is also the well known stochastic
calculus approach developed by Meyer \cite{meyer76} 
in which one works  with a general
semimartingale  $\{X_t\,:\, t\geq 0\}$, and defines   
$\Lambda=\{\lambda^{x}_t\,:\,x\in\R ,\, t\geq
0\}$ with respect to the Lebesgue measure from the formula
\begin{equation}
\label{eq2}
\int_0^t f(X_s)\,d<\hskip-1mm X^c\hskip-1mm>_s =
\int_{-\infty}^\infty f(x) \lambda^{x}_t\,dx.
\end{equation}
Of course, in the particular case when $d<\hskip-1mm X^c\hskip-1mm>_s=ds,$ i.e.,
$X^c$ is a Brownian motion, then the definitions of $L$ and $\lambda$ coincide.
In other cases, e.g., if $X^c\equiv 0,$ they will differ.

In this paper we focus on the potential
theoretic approach applicable in the Markovian case in which
the local times are defined as additive
functionals whose  $p$-potentials are equal to $p$-resolvent
kernels of $X.$  Local times can hereby be interpreted  as the
increasing processes in the Doob-Meyer decompositions of certain
submartingales. Considering the $p$-resolvent
kernels and passing to the limit, in an adequate manner, as $p\to
0, $ we obtain a formula (\ref{c1}), which clearly
extends Tanaka's original formula for the local times of Brownian
motion to those of the
symmetric  $\al$-stable processes, $\al\in(1,2],$
already obtained by T. Yamada \cite{yamada97} and further developed in
K. Yamada \cite{yamada02}. Our  approach may be simpler and may help to
make these results better known to probabilists working with L\'evy
processes.

The formula (\ref{c1}) below and its counterparts about decompositions of powers of symmetric
$\al$-stable L\'evy processes show at the same time similarities and
differences with the well known formulae for Brownian motion (see, in
particular, Chapter 10 in \cite{yor97} concerning the principal values
of Brownian local times). We hope that the Tanaka representation of
the local times in (\ref{c1}) may be useful to gain
some better understanding for the Ray-Knight theorems of the local
times of $X$ as presented in Eisenbaum et al. \cite{eisenbaumzhi00},
since in the Brownian case, Tanaka's formula has been such a powerful
tool for this purpose, see, e.g., Jeulin \cite{jeulin84}.

We now state the main formulae and results for the symmetric
$\al$-stable L\'evy process $X=\{X_t\}$. To be precise, we take
$X$ to satisfy 
$$
\E\left(\exp({\rm i}\lambda X_t)\right)=\exp(-t|\lambda|^\al),\quad \lambda\in\R,
$$
in particular, for $\al=2,$ $X$ equals $\sqrt{2}$ times a standard BM.
General criteria can be applied to
verify that $X$ possesses a jointly continuous family of
  local times $\{L^x_t\}$ satisfying (\ref{eq1}).
The constants $c_i$ 
appearing
below and  later in the paper will be computed
precisely in Section 5; clearly, they depend on the index $\al$ and/or
the exponent $\gamma$.

\vskip.5cm

\noindent
{\bf 1)} {\sl For all $t\geq 0$ and $x\in\R$
\begin{equation}
\label{c1}
|X_t-x|^{\alpha-1}= |x|^{\alpha-1} + N^x_t+ c_1\,L^x_t,
\end{equation}
where $N^x$ is a martingale 
such that for $0\leq \ga<\alpha/(\alpha-1),$
especially for $\gamma=2,$
\begin{equation}
\E\left(\sup_{s\leq t}|N^x_s|^\ga\right)<\infty.
\end{equation}
Moreover, the continuous
increasing process associated with
$N^x $ is 
\begin{equation}
\label{c2}
<\hskip-1mm N^x\hskip-1mm >_t\,:= c_2\,\int_0^t\frac{ds}{|X_s-x|^{2-\alpha}}.
\end{equation}
}
\vskip.5cm

\noindent
{\bf 2)} {\sl For $\alpha-1<\gamma<\alpha$  the submartingale
$\{|X_t-x|^\ga\}$ has the decomposition 
\begin{equation}
\label{e25}
|X_t-x|^\ga = |x|^\ga+N^{(\ga)}_t + A^{(\ga)}_t,
\end{equation}
where $N^{(\ga)} $ is a martingale and $A^{(\ga)}$ is the increasing process
given by
\begin{equation}
\label{e251}
A^{(\ga)}_t\,:=c_3\,\int_0^t \frac{ds}{|X_s-x|^{\alpha-\ga}}.
\end{equation}
}
\vskip.5cm

\noindent
{\bf 3)} {\sl  For $0<\gamma<\alpha-1$
the process $\{|X_t-x|^\ga\}$ is not a semimartingale but for $(\alpha-1)/2<\gamma<\alpha-1$ it 
is a Dirichlet process with the canonical decomposition
\begin{equation}
\label{e26}
|X_t-x|^\ga = |x|^\ga+N^{(\ga)}_t + A^{(\ga)}_t,
\end{equation}
where $N^{(\ga)} $ is a martingale and $A^{(\ga)},$ which has zero 
quadratic variation, is given by the
principal value integral
\begin{equation}
\label{c4}
A^{(\ga)}_t\,:=c_4\ {\rm p.v.}\,\int_0^t \frac{ds}{|X_s-x|^{\alpha-\ga}}\,:=
c_4\,\int\,\frac{dz}{|z|^{\alpha-\ga}}\left(L^{x+z}_t-L^x_t\right).
\end{equation}
}

The paper is organized so that in Section 2 some preliminaries about
symmetric L\'evy processes including their
generators and some variants of the It\^o formula are presented. In
Section 3
we derive the Tanaka formula for general 
symmetric L\'evy processes admitting local times. The above stated
results for symmetric stable L\'evy processes are proved and extended
in Section 4. In Section 5 we compute explicitly the
constants $c_i$ featured above and also further ones appearing especially
in Section 4. This is done by exhibiting some close relations between
these constants and the known expressions of the moments
$\E(|X_1|^\ga)$ where $X_1$ denotes a standard symmetric $\al$-stable variable.
In Section 6, we consider, instead of $|X_t-x|^{\ga}$, the process
$\{(X_t-x)^{\ga,*}\},$ where 
$$
a^{\ga,*}:={\rm sgn}(a)\,|a|^\ga,
$$
is the symmetric power of order $\ga,$ and we determine the parameter
values for which these
processes are semimartingales or Dirichlet processes, thus completing
results 1), 2) and 3) above.

\section{Preliminaries on symmetric L\'evy processes}
\label{sec2}

Throughout this paper, we consider a real-valued symmetric L\'evy
process $X=\{X_t\}$ and, if nothing else is stated, we assume $X_0=0.$ 
The L\'evy exponent $\Psi$ of $X$  is a non-negative
symmetric function such that
\begin{equation}
\label{e2001}
\E\left(\exp({\rm i}\xi X_t)\right)=\E\left(\cos(\xi X_t)\right)=
\exp\left(-t\Psi(\xi)\right).
\end{equation}
The L\'evy measure $\nu$ of $X $ satisfies, as is well known,
the integrability condition
$$
\int_{-\infty}^{\infty}\,(1\wedge z^2)\,\nu(dz)<\infty.
$$
By symmetry,
$\nu(A)=\nu(-A)$ for any $A\in\cB,$ the Borel $\sigma$-field  on $\R;$ hence,
\begin{eqnarray}
\label{e201}
&&\nonumber
\Psi(\xi)=\frac 12\,\sigma^2\,\xi^2-
\int_{-\infty}^{\infty}\left({\rm e}^{{\rm i}\,\xi\,z}-1-
{\rm i}\,\xi\,z\,{\bf 1}_{\{|z|\leq 1\}}\right)\,
\nu(dz)\\
&&\hskip1cm
= \frac 12\,\sigma^2\,\xi^2+2
\int_{0}^{\infty}\left(1-\cos(\xi z)\right)
\,
\nu(dz).
\end{eqnarray}
Recall also (see, e.g., Ikeda and Watanabe \cite{ikedawatanabe81}
p. 65) that $X$ admits the Brownian-Poisson representation
\begin{equation}
\label{e2010}
X_t=\sigma\, B_t+\int_{(0,t]}\int_{\{|z|\geq 1\}}z\,\Pi(ds,dz)+
\int_{(0,t]}\int_{\{|z|<1\}}z\,(\Pi-\pi)(ds,dz),
\end{equation}
where the Brownian motion $B$ and the Poisson random measure  $\Pi$  with
the intensity
$$
\pi(ds,dz):=\E(\Pi(ds,dz))=ds\,\nu(dz)
$$
are independent. Due to the symmetry of $\nu,$ 
the generator of $X$ can be written as 
\begin{eqnarray}
\label{gene1}
&&
\nonumber
\cG f(x):=\cG^B f(x)+\cG^\Pi f(x)\\
&&
\nonumber
\hskip1.2cm
:=
\frac 12\, \sigma^2\, f''(x)+
\int_\R \left(f(x+y)-f(x)-f'(x)\,y\,{\bf
  1}_{\{|y|<1\}}\right)\,\nu(dy)\\
&&
\hskip1.2cm
=
\frac 12\, \sigma^2\, f''(x)+
\int_\R \left(f(x+y)-f(x)-f'(x)\,y\right)\,\nu(dy).
\end{eqnarray}
where $\cG$ acts on regular functions $f$ in particular those in 
the Schwartz space $S(\R)$ of rapidly
decreasing functions. Given a smooth function $f,$ 
the predictable form of the It\^o formula (see Ikeda and Watanabe
\cite{ikedawatanabe81} and K. Yamada \cite{yamada02}) writes 
\begin{eqnarray}
\label{ito}
&&\hskip-1cm
f(X_t)-f(X_0)-\int_0^t\cG f(X_s)\, ds \\
&&
\nonumber
=
\sigma\,\int_0^tf'(X_s)\, dB_s
+ \int_0^t\int_\R\left(f(X_{s-}+z)-f(X_{s-})\right)\,
(\Pi-\pi)(ds,dz).
\end{eqnarray}
The formula (\ref{ito}) connects with the It\^o formula
for semimartingales, as developed by Meyer \cite{meyer76}, 
and displayed as
\begin{eqnarray}
\label{meyer1}
&&\hskip-1cm
\nonumber
f(X_t)=f(X_0)+\int_0^tf'(X_{s-})\, dX_s+\frac {\sigma^2}2 \int_0^tf''(X_{s})\, ds\\&&
\hskip2cm
+\sum_{0<s\leq t}
\left(f(X_s) -f(X_{s-}) - f'(X_{s-})\Delta X_s\right).
\end{eqnarray}
The sum of jumps $\displaystyle{\sum_{0<s\leq t}(\dots)}$ may be
compensated by
$\displaystyle{\int_0^t\cG^\Pi f(X_s)\, ds},$ and, hence,  we have 
recovered the integrated form of (\ref{gene1}):
$$
\int_0^t\cG f(X_s)\, ds=\int_0^t\cG^B f(X_s)\, ds+ \int_0^t\cG^\Pi f(X_s)\, ds.
$$
We record also a more general compensator formula employed later in the
paper. For this, let
$\Phi:\R\times\R\mapsto \R_+$ be a Borel measurable function. Then  
\begin{eqnarray}
\label{com1}
&&\nonumber
\E\left(\sum_{0<s\leq t}
\Phi(X_{s-},X_{s}){\bf 1}_{\{\Delta X_s\not=0\}}\right)\\
&&
\hskip2cm
=
\E\left(\int_0^t\,\int_{\R\setminus\{0\}}\pi(ds,dz)
\Phi(X_{s},X_{s}+z)\right).
\end{eqnarray}


\section{Local times for symmetric L\'evy processes}
\label{sec3}

From now on, we assume that 
\begin{equation}
\label{e202}
\int_{-\infty}^{\infty}\frac{1}{1+\Psi(\xi)}\,d\xi<\infty.
\end{equation}
From standard Fourier arguments (see Bertoin \cite{bertoin96} and, e.g.,
Borodin and Ibragimov \cite{borodinibragimov95} p. 67) one can show
the existence of  a jointly measurable family of local
times $\{L^x_t\,:\, x\in\R,t\geq 0\}$ satisfying 
for every Borel-measurable function $f:\R\mapsto\R_+$
the occupation time
formula
$$
\int_0^t ds\, f(X_s)=
\int_{-\infty}^\infty f(x) L^x_t\,dx.
$$
For the condition (expressed in terms of the function $v$ in
(\ref{e240})) 
under which $(t,x)\mapsto L^x_t$ is continuous, see
Bertoin \cite{bertoin96} p. 148. In particular, the condition holds for
symmetric $\al$-stable L\'evy processes; in fact it was shown by
Boylan \cite{boylan64}, see also Getoor and Kesten
\cite{getoorkesten72}, that
\begin{equation}
\label{f30}
|L^{x+y}_t-L^x_t|\leq K_t\, |y|^\theta
\end{equation}
for any $\theta<(\alpha-1)/2$ and some random constant $K_t.$

Our approach toward a Tanaka formula for these local times is based on
the potential theoretic construction which we now develop. It is well
known, see Bertoin \cite{bertoin96} p. 67,  that for any $p>0$ 
\begin{equation}
\label{e21}
u^{(p)}(x)=\frac 1\pi\int_0^\infty\frac{\cos(\xi
  x)}{p+\Psi(\xi)}
\,d\xi
\end{equation}
is a continuous version of the density of the resolvent
$$
U^{(p)}(0,dx)=\E_0\left(\int_0^\infty {\rm e}^{-p\,t}{\bf 1}_{\{X_t\in
  dx\}}\, dt\right).
$$
Moreover, for every $x$  the local time $\{L^x_t\}$ can be chosen
  as a continuous additive functional such that
\begin{equation}
\label{eq3}
u^{(p)}(y-x)=\E_y\left(\int_0^\infty {\rm e}^{\,-p\, t}
d_t L^x_t\right).
\end{equation}
From (\ref{eq3}) we deduce the Doob-Meyer decomposition given in the next

\begin{proposition}
\label{le1}
For every fixed $x$
\begin{equation}
\label{eq4}
u^{(p)}(X_t-x)=u^{(p)}(X_0-x)
+M^{(p,x)}_t+p\,\int_0^t u^{(p)}(X_s-x)ds -L^{x}_t,
\end{equation}
where $M^{(p,x)}$ is a martingale with respect
to the natural filtration $\{\cF_t\}$ of $X.$
Moreover, for every fixed $t,$ both the martingale $\{M^{(p,x)}_s\,:\,
s\leq t\}$ and the random variable $L^{x}_t$ belong to BMO; in particular,
$L^{x}_t$ has some exponential moments.
\end{proposition}
\begin{proof}
Straightforward computations using the Markov property show that
for $y=X_0$
\begin{eqnarray*}
&&
\E_y\left(\int_0^\infty {\rm e}^{\,-p\, t} d_t L^{x}_t\,|\, \cF_s\right)
=
\int_0^s{\rm e}^{\,-p\, t} d_t L^{x}_t
+{\rm e}^{\,-p\, s}
u^{(p)}(X_s-x),
\end{eqnarray*}
which together with an integration by parts yields
(\ref{eq4}). We leave the proofs of the remaining assertions to the reader.
\end{proof}
A variant of the Tanaka formula shall now be obtained by letting $p\to
0$ in (\ref{eq4}). The result is stated in Proposition \ref{prop11} but first we need
an important ingredient.

\begin{lemma}
\label{lemma1}
For every $x\in\R$
\begin{equation}
\label{e240}
\lim_{p\to 0}\left( u^{(p)}(0)-u^{(p)}(x)\right)
= \frac 1\pi\,\int_0^\infty \frac{1-\cos(\xi x)}{\Psi(\xi)}\, d\xi
=:v(x).
\end{equation}
\end{lemma}

\begin{proof}
The statement follows from (\ref{e21}) by dominated convergence because
(cf. (\ref{e201}))
$$
\int^\infty_1 \frac{1}{\Psi(\xi)}\,d\xi<\infty,\quad {\rm and}\quad
\int_0^1 \frac{\xi^2}{\Psi(\xi)}\,d\xi<\infty.
$$
Notice also that $v$ is continuous. 
\end{proof}

The formula (\ref{e24}) below generalizes in a sense the Tanaka formula
for Brownian motion to symmetric L\'evy processes.
In the next section we study the particular case of symmetric stable processes.

\begin{proposition}
\label{prop11}
Let $v$ be the function introduced
in (\ref{e240}) and $M^{(p,x)}$ the martingale defined in
Proposition \ref{le1}. Then
\begin{equation}
\label{e24}
v(X_t-x)=v(x)+
\tilde N^{x}_t+
L^x_t,
\end{equation}
where $\tilde N^x_t:=-\lim_{p\to 0}M^{(p,x)}_t$ defines a martingale.

\end{proposition}
\begin{remark}
{\rm 
Standard results about martingale additive functionals of $X$ yield
the following representations 
\begin{eqnarray*}
&&
\tilde N^x_t
={\sigma} \int_0^t v'(X_s-x)\, dB_s
\\
&&
\hskip2cm
+\int_{(0,t]}\int_{\R}
\left( v(X_{s-}-x+z)-v(X_{s-}-x)\right)\,
\,(\Pi-\pi)(ds,dz),
\end{eqnarray*}
and
\begin{eqnarray*}
&&
<\tilde N^x>_t=\sigma^2\,\int_0^t (v'(X_s-x))^2\, ds\\
&&
\hskip3cm
 +
\int_0^t\int_{\R}
\left( v(X_{s}-x+z)-v(X_{s}-x)\right)^2\,
\,\pi(ds,\,dz),
\end{eqnarray*}
where $v'$ is a weak derivative of $v.$
}
\end{remark}

\begin{proof}
Consider the identity (\ref{eq4}). 
Let therein $p\to 0$ and use Lemma \ref{lemma1} to obtain
\begin{equation}
\label{e22}
v(X_t-x)=v(x)-
\lim_{p\to 0}\left(M^{(p,x)}_t+
p\int_0^t u^{(p)}(X_s-x)ds\right) +L^x_t.
\end{equation}
From (\ref{e21})
$
u^{(p)}(y)\leq u^{(p)}(0),
$
and, consequently,
\begin{equation}
\label{e221}
0\leq p\int_0^t u^{(p)}(X_s-x)ds
\leq p\,u^{(p)}(0)\,t.
\end{equation}
Next we show that
\begin{equation}
\label{e231}
\lim_{p\to 0}p\,u^{(p)}(0)= 0.
\end{equation}
Indeed, using (\ref{e21}) again,
\begin{eqnarray}
\label{e2301}
\nonumber
&&
p\,u^{(p)}(0)=\frac 1\pi\ \int_0^\infty \frac{p\,d\xi}{p+\Psi(\xi)}\\
&&\hskip1.6cm
\leq\frac 1\pi\ \int_0^1 \frac{p\,d\xi}{p+\Psi(\xi)}
+\frac p\pi\ \int_1^\infty \frac{d\xi}{\Psi(\xi)},
\end{eqnarray}
and (\ref{e231}) results by  dominated convergence.
Hence, (\ref{e22}) yields (\ref{e24}) with
$\tilde N^x$ as claimed. It remains to prove that $\tilde N^x$ is a martingale.
For this it is enough to show that 
\begin{equation}
\label{e2311}
\E\left(|\tilde N^x_t-M^{(p,x)}_t|\right)\to 0\quad {\rm as}\ p\to 0.
\end{equation}
To prove (\ref{e2311}) consider 
\begin{eqnarray*}
&&\hskip-1cm
|\tilde N^x_t-M^{(p,x)}_t|
\leq 
p\int_0^t u^{(p)}(X_s-x)ds
+|v(x)-(u^{(p)}(0)-u^{(p)}(x))|\\
&&
\hskip3cm
+|v(X_t-x)-(u^{(p)}(0)-u^{(p)}(X_t-x))|.
\end{eqnarray*}
From (\ref{e221}) and (\ref{e231}), the integral term goes to 0 as
$p\to 0.$
Next, by Fubini's theorem  and (\ref{e2001})
\begin{eqnarray*}
&&
\E\left(\left|v(X_t-x)-(u^{(p)}(0)-u^{(p)}(X_t-x))\right|\right)\\
&&\hskip2cm
=\frac{1}{\pi}\,
p\int_0^\infty\frac{1-\E(\cos(\xi(X_t-x)))}{\Psi(\xi)(p+\Psi(\xi))}\,
d\xi\\
&&\hskip2cm
=\frac{1}{\pi}\,
p\int_0^\infty\frac{1-\cos(\xi\,x)\,\exp(-t\Psi(\xi))}{\Psi(\xi)(p+\Psi(\xi))}\,
d\xi\\
&&\hskip2cm
\leq\frac{1}{\pi}\,
p\left(\int_0^\infty\frac{1-\cos(\xi\,x)}
{\Psi(\xi)(p+\Psi(\xi))}\,d\xi
+
\int_0^\infty\frac{t\,\Psi(\xi)}
{\Psi(\xi)(p+\Psi(\xi))}\,
d\xi\right).
\end{eqnarray*}
Applying the dominated convergence theorem for the first term above and
(\ref{e2301}) for the second one give
$$
\lim_{p\to 0}
\E\left(\left|v(X_t-x)-(u^{(p)}(0)-u^{(p)}(X_t-x))\right|\right)
=0,
$$
completing the proof.
\end{proof}

\begin{example}
{\rm
For standard Brownian motion $B$  we have
$$
u^{(p)}(x)=\frac{1}{\sqrt{2p}}\ {\rm e}^{\,-\sqrt{2p}|x|}.
$$
Consequently,
$$
v(x):=\lim_{p\to 0}\left( u^{(p)}(0)-u^{(p)}(x)\right)=|x|.
$$
and the formula (\ref{e24}) takes the familiar form
$$
|B_t-x|=|x|+N^{x}_t+L^x_t
$$
where
\begin{eqnarray*}
&&N^{x}_t=\lim_{p\to 0} \int_0^t {\rm e}^{\,-\sqrt{2\,p}|B_s-x|}\,
{\rm sgn}(B_s-x)\,dB_s\\
&&\hskip.7cm
=\int_0^t {\rm sgn}(B_s-x)\,dB_s.
\end{eqnarray*}
}
\end{example}

\section{Symmetric $\al$-stable L\'evy processes}
\label{sec4}

Let $X=\{X_t\},\, X_0=0,$ denote the symmetric $\alpha$-stable
process with the L\'evy exponent
$$
\Psi(\xi)=|\xi|^{\alpha}, \quad \alpha\in(1,2).
$$
We remark that the condition (\ref{e202}) is satisfied, and also that 
the local time of $X$ has a jointly continuous version, as is
discussed in Section 3.
For clarity,
we have excluded the Brownian motion 
from our study. However, the corresponding results for Brownian
motion may be recovered by letting $\al\to 2.$
Recall also that $\E(|X_t|^\ga)<\infty$ for $\ga< \al,$ and that 
the L\'evy measure is
\begin{equation}
\label{c5}
\nu(dz)= c_5(\al)\, |z|^{-\alpha-1}\, dz, \quad \alpha\in(1,2).
\end{equation}
The function $v$ introduced in Lemma \ref{lemma1}
is in the present case given by
\begin{equation}
\label{c6}
v(x)= c_6(\al)\,|x|^{\alpha-1}.
\end{equation}

The results announced in the Introduction are now presented again and proven in a more complete form through the following
three propositions. The first one treats the claim 1) in the Introduction.

\begin{proposition}
\label{prop1} {\bf a)} For fixed $x$
\begin{equation}
\label{e222}
c_6(\alpha)\,\left( |X_t-x|^{\alpha-1} -|x|^{\alpha-1}\right)= \tilde N^x_t+ L^x_t\end{equation}
where $\{\tilde N^x_t\}$
is a square integrable martingale. In fact, for all
$0\leq \ga<\alpha/(\alpha-1),$
especially for $\gamma=2,$
\begin{equation}
\label{e2221}
\E\left(\sup_{s\leq t}|\tilde N^x_s|^\ga\right)<\infty.
\end{equation}
Moreover, the continuous
increasing process associated with 
$\tilde N^x $ is 
\begin{equation}
\label{c7}
<\hskip-1mm \tilde N^x\hskip-1mm >_t:= c_7(\alpha)\,\int_0^t\frac{ds}{|X_s-x|^{2-\alpha}}.\end{equation}
{\bf b)} For every $t$ and $x$ the variable $L^x_t$ belongs to BMO;
in fact, for all $s\leq t$
\begin{equation}
\label{e2222}
\E\left(L^x_t-L^x_s\,|\,\cF_s\right)\leq K_{\alpha,t}
\end{equation}
for some constant $K_{\alpha,t}$ which does not depend on $s.$
\end{proposition}

\begin{proof} The fact that $\tilde N^x$ is a martingale
is clear from Proposition \ref{prop11}.
Because $L^x_t$ has some exponential moments (cf. Proposition \ref{le1}), it is seen easily from(\ref{e222}) that for $\ga>0$
$$
\E(|\tilde N^x_t|^\ga)<\infty
$$
if
$$
\E\left(|X_t-x|^{\ga(\alpha-1)}\right)<\infty,
$$
which is true for  $\ga(\alpha-1)<\alpha.$ Consequently, an extension of the
Doob-Kolmogorov inequality,
gives (\ref{e2221}). The martingale $\tilde N^x$ has no continuous
martingale part. Hence, letting
$$
[\tilde N^x]_t:=\sum_{s\leq t}(\Delta \tilde N^x_s   )^2
:=(c_6(\al))^2\sum_{s\leq
  t}\left(|X_s-x|^{\alpha-1}-|X_{s-}-x|^{\alpha-1}\right)^2
$$
it holds that $\{\,(\tilde N^x_t)^2-[\tilde N^x]_t\,\}$ is a martingale.  Consequently,  $<\hskip-1mm \tilde N^x\hskip-1mm>$
can be obtained as the dual predictable projection of $[\tilde N^x],$ and 
from the L\'evy system of $X$, e.g., (\ref{com1}),  we get
$$
<\hskip-1mm \tilde N^x\hskip-1mm >_t=
(c_6(\al))^2\,c_5(\al)\,\int_0^tds\int_\R \frac{dy}{|y|^{\alpha+1}}\left(|X_{s-}-x+y|^{\alpha-1}-|X_{s-}-x|^{\alpha-1}\right)^2.$$
Putting $z=X_{s-}-x$ and introducing $y=zu$ the latter integral takes
the form 
$$
\int_\R \frac{dy}{|y|^{\alpha+1}}\left(|z+y|^{\alpha-1}-
|z|^{\alpha-1}\right)^2
=
\frac{1}{|z|^{2-\alpha}}\int_\R \frac{du}{|u|^{\alpha+1}}\left(|1+u|^{\alpha-1}-
1\right)^2.
$$
Consequently, $<\hskip-1mm \tilde N^x\hskip-1mm >$ is as claimed. To prove the second part of the proposition, notice that by the martingale property
\begin{eqnarray*}
&&\E\left(L^x_t-L^x_s\,|\,\cF_s\right)
=c_6(\alpha)\ \E\left(|X_t-x|^{\alpha-1}
-|X_s-x|^{\alpha-1}\,|\,\cF_s\right)\,\\
&&
\hskip2.9cm
\leq c_6(\alpha)\
\E\left(|X_t-X_s|^{\alpha-1}\,|\,\cF_s\right)\\
&&
\hskip2.9cm
\leq c_6(\alpha)\
\E\left(|X_{t-s}|^{\alpha-1}\right)
\\
&&
\hskip2.9cm
\leq K'_\alpha\,t^{(\alpha-1)/\alpha},
\end{eqnarray*}
where also the scaling property and the inequality
$$
|x^p-y^p|\leq |x-y|^p,\quad 0<p\leq 1,
$$
are used.
\end{proof}

The following corollary plays the same r\^ole for $X$ as the
classical It\^o-Tanaka formula
plays for Brownian motion. In fact, a large part of this paper discusses
for which functions the identity (\ref{itotanaka}), or some variant of
it is valid.
\begin{corollary}
\label{cor1}
Let $f$ be a bounded Borel function with compact support and define
$$
F(y):=\int dx f(x)\,|y-x|^{\al-1}.
$$
Then
\begin{equation}
\label{itotanaka}
F(X_t)=F(0)+\int dx f(x)\, N^x_t+c_1(\al)\int_0^t ds\, f(X_s)
\end{equation}
expresses the canonical semimartingale decomposition of $\{F(X_t)\}$
with\break $\{\int dx f(x)\, N^x_t\}$ a martingale.
\end{corollary}
\begin{proof} It suffices to integrate both sides of 
(\ref{e222}) (or rather (\ref{c1})) with respect to the measure $f(x)\,dx.$
\end{proof}
\begin{remark} 
{\bf a)}  {\rm In K. Yamada \cite{yamada02} the representation (\ref{e222})
  of the local time (or Tanaka's formula for symmetric $\al$-stable
  processes) is derived using the so called ``mollifier'' approach 
as in Ikeda and Watanabe  \cite{ikedawatanabe81} in the Brownian
motion case. In this case the martingale is given by 
$$
\tilde N^x_t=c_6(\al) \,\int_{(0,t]}\int_{\R}
\left( |X_{s-}-x+z|^{\al-1}-|X_{s-}-x|^{\al-1}\right)\,
\,(\Pi-\pi)(ds,dz),
$$
where $\Pi$ and $\pi$ are the Poisson random measure and the
corresponding  intensity measure, respectively, associated with $X$.  
}
\hfill\break\hfill
{\bf b)}  {\rm
The inequality (\ref{e2222}) holds for all symmetric L\'evy
processes having local times. Indeed, it is proved in 
Bertoin \cite{bertoin96} p. 147 Corollary 14  that the function  $v$
defined in (\ref{e240}), Lemma \ref{lemma1}, induces a metric on
$\R,$ and, in particular, the triangle inequality holds. Consequently, 
\begin{eqnarray*}
&&\E\left(L^x_t-L^x_s\,|\,\cF_s\right)
\leq  \E(v(X_t-X_s))
=  \E(v(X_{t-s}))
\leq  \E(v(X_{t}))
<\infty
\end{eqnarray*}
because
\begin{eqnarray*}
&&
\E(v(X_t))=\frac 1\pi\,\int_0^\infty \frac{1-
\exp(-t\Psi(\xi))}
{\Psi(\xi)}\, d\xi\\
&&
\hskip1.8cm
\leq \frac 1\pi\,\left( 
t
+ \int_1^\infty \frac{d\xi}{\Psi(\xi)}\right)<\infty.
\end{eqnarray*} 
}
\hfill\break\hfill
{\bf c)}  {\rm We leave it to the reader to establish a 
  version of Corollary \ref{cor1} for general symmetric L\'evy processes. 
}
\end{remark}

\begin{proposition}
\label{prop2}
For a given $x$ and $\alpha-1<\gamma<\alpha$  the submartingale
$\{|X_t-x|^\ga\,:\, t\geq 0\}$ has the decomposition
\begin{equation}
\label{e25}
|X_t-x|^\ga = |x|^\ga+N^{(\ga)}_t + A^{(\ga)}_t,
\end{equation}
where $N^{(\ga)} $ is a martingale and $A^{(\ga)}$ is the increasing process
given by
\begin{equation}
\label{c3}
A^{(\ga)}_t=c_3(\al,\ga)\,\int_0^t \frac{ds}{|X_s-x|^{\al-\ga}}.
\end{equation}
Moreover, when $\alpha-1\leq \ga\leq \alpha/2$ the increasing
process $<\hskip-1mm N^{(\ga)}\hskip-1mm >$ is of the form
\begin{equation}
\label{c8}
<\hskip-1mm N^{(\ga)}\hskip-1mm >_t=c_8(\al,\ga)\,\int_0^t \frac{ds}{|X_s-x|^{\al-2\ga}}.\end{equation}
\end{proposition}
\begin{proof}
Formula (\ref{e25}) is obtained by integrating both sides of equation (\ref{e222})
(or (\ref{c1}) taken at level $z$ with respect to the measure
$dz/|z-x|^{\al-\ga}.$ The form of the left hand side is obtained from
the scaling argument.
Because $A^{(\ga)}$ is continuous the computation for finding
$<\hskip-1mm N^{(\ga)}\hskip-1mm >$ is very similar to
the computation of $<\hskip-1mm \tilde N^x\hskip-1mm >$ in the proof of
Proposition  \ref{prop1}. We have
\begin{equation}
\label{e260}
<\hskip-1mm N^{(\ga)}\hskip-1mm >_t=
\int_0^tds\int_\R \nu(dy)\left(|X_{s-}-x+y|^{\ga}-
|X_{s-}-x|^{\ga}\right)^2,
\end{equation}
which easily yields  (\ref{c8}).
\end{proof}

For the next proposition, we recall the notion of Dirichlet
process, that is a process which can decomposed uniquely as the sum of a 
local martingale and a continuous process with zero quadratic
variation (see, e.g., F\"ollmer \cite{follmer81}, Fukushima \cite{fukushima80}).

\begin{proposition}
\label{prop3}
{\bf a)}  For $0<\gamma<\alpha-1$
the process $|X-x|^\ga$
is not a semimartingale.  \hfill\break
{\bf b)} For $(\alpha-1)/2<\gamma<\alpha-1$
the process $|X-x|^\ga$
is a Dirichlet process with the canonical decomposition
\begin{equation}
\label{e261}
|X_t-x|^\ga = |x|^\ga+N^{(\ga)}_t + A^{(\ga)}_t,
\end{equation}
where $N^{(\ga)} $ is a martingale and $A^{(\ga)}$ is given by the
principal value integral
\begin{eqnarray}
\label{ec4}
&&
\nonumber
A^{(\ga)}_t=c_4(\al,\ga)\ {\rm p.v.}\, \int_0^t
\frac{ds}{|X_s-x|^{\al-\ga}}\\
&&\hskip.9cm
=
c_4(\al,\ga)\,\int_\R\,\frac{dz}{|z|^{\al-\ga}}\left(L^{x+z}_t-L^x_t\right).
\end{eqnarray}
Moreover, the increasing process $<\hskip-1mm N^{(\ga)}\hskip-1mm >$
is as given in (\ref{c8}).
\end{proposition}

\begin{proof}  a) We take $x=0$ and adapt the argument in Yor \cite{yor78}
  applied therein for continuous martingales.  Assume that
  $Y_t:=|X_t|^\ga,\, \ga<\alpha-1,$ defines a
  semimartingale. Then
$$
|X_t|^{\alpha-1}=Y_t^{\,\theta}
$$
with $\theta= \ga/(\alpha-1)>1,$ and It\^o's formula for semimartingales
(notice that $Y^c\equiv 0$) gives
\begin{equation}
\label{e27}
Y_t^\theta=\int_0^t\theta\,Y_{s-}^{\theta-1}\,dY_s+\Sigma_t,
\end{equation}
where
$$
\Sigma_t:=
\sum_{0<s\leq t}
\left(Y_s^\theta-Y_{s-}^\theta-\theta\,Y_{s-}^{\theta-1}\Delta
Y_s\right).
$$
The argument of the proof is that under the above assumption
the local time
\begin{equation}
\label{e2611}
L^{0}_t\equiv\int_0^t {\bf 1}_{\{X_{s-}=0\}}\, d|X_s|^{\alpha-1}
\end{equation}
would be equal to zero. To derive this  contradiction notice from
(\ref{e27}) and (\ref{e2611}) that
\begin{eqnarray*}
&&
L^{0}_t=\int_0^t {\bf 1}_{\{Y_{s-}=0\}}\, dY_s^{\theta}
=\int_0^t{\bf 1}_{\{Y_{s-}=0\}}\,d\Sigma_s.
\end{eqnarray*}
But because  $\Sigma$ is a purely discontinuous increasing process and $L^0$ is
continuous this is possible only if $L^0\equiv 0,$
which cannot be the case; thus proving that $Y$ is not a semimartingale.
\hfill\break\hfill
b) To prove (\ref{e261}) we consider formula (\ref{c1}) at levels $x+z$ and $x$ and write
\begin{eqnarray}
\label{f26}
&&
\int_\R\frac
{dz}{|z|^{\al-\ga}}\left(|X_t-(x+z)|^{\al-1}-|X_t-x|^{\al-1}\right)\\
&&
\nonumber
\hskip2cm
=
\int_\R\,\frac{dz}{|z|^{\al-\ga}}\left(N^{x+z}_t-N^x_t\right)+
c_1(\al)\,\int_\R\,\frac{dz}{|z|^{\al-\ga}}\left(L^{x+z}_t-L^x_t\right).
\end{eqnarray}
The integral on the left hand side is well defined since by scaling 
$$
\int_\R\frac
{dz}{|z|^{\al-\ga}}\left(|X_t-(x+z)|^{\al-1}-|X_t-x|^{\al-1}\right)
=
|X_t-x|^{\ga}\,r(\al,\ga)
$$
with
$$
r(\al,\ga):=\int_\R\frac
{dz}{|z|^{\al-\ga}}\left(|1-z|^{\al-1}-1\right),
$$
which is an absolutely convergent integral. Next notice that the
principal value integral on the right hand side of (\ref{f26}) is
well defined by the H\"older continuity in $x$ of the local times 
(cf. (\ref{f30})). It also follows that the first integral on the
right hand side of (\ref{f26}) is meaningful and, by Fubini's
theorem, it is a martingale. In  Fitzsimmons and Getoor
\cite{fitzsimmonsgetoor92a} it is proved that $A^{(\ga)}$ 
has zero $p$-variation for 
$p>p_o:=(\al-1)/(\al-1-\gamma)$ and infinite  $p$-variation for $p<p_o.$ 
Hence the claimed Dirichlet process decomposition follows with 
\begin{equation}
\label{e270}
c_4(\al,\ga)=c_1(\al)/r(\al,\ga).
\end{equation}
\end{proof}

\section{Explicit values of the constants}
\label{sec5}

An important ingredient in the computation of the explicit values of the
constants is the formula for absolute moments of symmetric
$\al$-stable, $\al\in(1,2),$ random variables due to Shanbhag and Sreehari
\cite{shanbhagsreehari77} (see also Sato \cite{sato99} p. 163,
Chaumont and Yor \cite{chaumontyor03} p. 110). To discuss
this briefly let 
\begin{description}
\item{$\circ$}\quad $Z$ be an exponentially distributed r.v. with mean 1,
\item{$\circ$}\quad $U$ a normally distributed r.v. with mean 0 and variance 1,
\item{$\circ$}\quad $X^{(\al)}$ a symmetric $\al$-stable r.v. with
  characteristic function $\exp(-|\xi|^\al),$
\item{$\circ$}\quad $Y^{(\al/2)}$ a positive $\al/2$-stable r.v. with
  Laplace transform $\exp(-\xi^{\al/2}).$
\end{description}

Assume also that these variables are independent. Then it is easily
checked that 
\begin{equation}
\label{e510}
\left(Z/Y^{(\al/2)}\right)^{\al/2}\ \rr\ Z 
\end{equation}
and 
\begin{equation}
\label{e520}
X^{(\al)}\ \rr\ \sqrt 2\, U\, \left(Y^{(\al/2)}\right)^{1/2}.
\end{equation}
From (\ref{e510}) we obtain for $\ga<\al/2$
$$
\E\left(\left(Y^{(\al/2)}\right)^{\ga}\right)=\frac{\Gamma(1-\frac{2\ga}{\al})}
 {\Gamma(1-\ga)},
$$
and, further, from (\ref{e520}) for $-1<\ga<\al$
\begin{equation}
\label{e511}
m_\ga:=\E\left(|X^{(\al)}|^\ga\right)=2^\ga\,\Gamma(\frac{1+\ga}{2})\,
\Gamma(\frac{\al-\ga}{\al})/\left(\sqrt{\pi}\, \Gamma(\frac{2-\ga}{2})\right).
\end{equation}

The constants with the associated reference numbers of the formulae where
they appear in the paper are summarized in the following table. 

$$
\vbox{\offinterlineskip\halign{&\vrule#&\strut\quad\hfil
$#$\quad\hfil\cr
\noalign{\hrule}
height10pt&\omit&&\omit&&\omit&\cr
&{\rm Constant}&&{\rm Value}&&{\rm Ref.}&\cr
height10pt&\omit&&\omit&&\omit&\cr
\noalign{\hrule height0.1pt}
height10pt&\omit&&\omit&&\omit&\cr
&
c_1(\al) 
&&
\left((\al-1)\,\pi\, m_{\al-1}\right)/\Gamma(1/\al)
&&
(\ref{c1}), (\ref{e61}), (\ref{e62})
&
\cr
height10pt&\omit&&\omit&&\omit&\cr
\noalign{\hrule height0.1pt}
height10pt&\omit&&\omit&&\omit&\cr
&
c_2(\al) 
&&
(2(\al-1)\, m_{2(\al-1)})/(\al\,m_{\al-2})
&&
(\ref{c2})
&
\cr
height10pt&\omit&&\omit&&\omit&\cr
\noalign{\hrule height0.1pt}
height10pt&\omit&&\omit&&\omit&\cr
&
c_3(\al,\ga) 
&&
(\ga\, m_\ga)/(\al\,m_{\ga-\al})
&&
(\ref{e251}), (\ref{c3})
&
\cr
height10pt&\omit&&\omit&&\omit&\cr
\noalign{\hrule height0.1pt}
height10pt&\omit&&\omit&&\omit&\cr
&
c_4(\al,\ga) 
&&
c_1(\al)/r(\al,\ga)
&&
(\ref{ec4}), (\ref{e270})
&
\cr
height10pt&\omit&&\omit&&\omit&\cr
\noalign{\hrule height0.1pt}
height10pt&\omit&&\omit&&\omit&\cr
&
c_5(\al) 
&&
\al/(2\,\Gamma(1-\alpha)\,\cos(\alpha\pi/2))
&&
(\ref{c5})
&
\cr
height10pt&\omit&&\omit&&\omit&\cr
\noalign{\hrule height0.1pt}
height10pt&\omit&&\omit&&\omit&\cr
&
c_6(\al) 
&&
(c_1(\al))^{-1}=
\left(2\pi c_5(\alpha-1)\right)^{-1}
&&
(\ref{c6})
&
\cr
height10pt&\omit&&\omit&&\omit&\cr
\noalign{\hrule height0.1pt}
height10pt&\omit&&\omit&&\omit&\cr
&
c_7(\al)
&&
c_2(\al)\,(c_6(\al))^{\,2}
&&
(\ref{c7})
&
\cr
height10pt&\omit&&\omit&&\omit&\cr
\noalign{\hrule height0.1pt}
height10pt&\omit&&\omit&&\omit&\cr
&
c_8(\al,\ga) 
&&
c_3(\al,2\ga)-2c_3(\al,\ga)
&&
(\ref{c8})
&
\cr
height10pt&\omit&&\omit&&\omit&\cr
\noalign{\hrule}}}
$$

We consider first the constant $c_3(\al,\ga)$ and, for clarity, recall  formula (\ref{e25}):
\begin{equation}
\label{e52}
|X_t-x|^\ga = |x|^\ga+N^{(\ga)}_t + A^{(\ga)}_t,
\end{equation}
with $\al-1<\ga<\al$ and 
$$
A^{(\ga)}_t=c_3(\al,\ga)\,\int_0^t \frac{ds}{|X_s-x|^{\al-\ga}}.
$$
Notice that letting  $\ga\downarrow \al-1$ yields, in a sense, 
\begin{equation}
\label{e53}
A^{(\al-1)}_t=c_1(\al)\,L^x_t,
\end{equation}
although, using the value in the table, $c_3(\al,\ga)\to 0.$ 
From  (\ref{e52}) it is seen that 
$f(y)=|y-x|^\ga$ belongs to the domain of the extended generator
$\cG,$ and, by scaling we obtain the following integral representation
$$
c_3(\al,\ga)=\int_\R\nu(dy)\left(|1+y|^\ga -1-\ga y\right).
$$
On the other hand, taking $x=0$ in (\ref{e52}), and using scaling
again together with  (\ref{e511}), we get 
$$
\E\left(|X_t|^\ga\right)
=c_3(\al,\ga)\int_0^tds\,\E\left(|X_s|^{\ga-\al}\right),  
$$
which is equivalent with
$$
t^{\ga/\al}\, m_\ga=c_3(\al,\ga)\,\frac{\al\, t^{\ga/\al}}{\ga}\,
m_{\ga-\al}
$$
hence,
$$
c_3(\al,\ga)=\ga\, m_\ga/\al\,m_{\ga-\al}.
$$ 
A similar argument leads to an expression for $c_1(\al).$ From
(\ref{e53}) we get
\begin{equation}
\label{e54}
\E\left(|X_t|^{\al-1}\right) =c_1(\al)\,\E\left(L^0_t\right).  
\end{equation}
We derive from (\ref{e54}) the existence of a constant $c_0(\al)$ such that
$$
\E\left(d_t L^0_t\right)=c_0(\al) dt \, t^{-1/\al},
$$
and it follows from (\ref{e53}) that
\begin{equation}
\label{e55}
m_{\al-1}=\al\,c_1(\al)\,c_0(\al)/(\al-1).
\end{equation}
We now compute $c_0(\al)$ to obtain $c_1(\al)$ from (\ref{e55}). For this
consider the identity (\ref{eq3}) for $x=y=0$
$$
u^{(p)}(0)=\E_0\left(\int_0^\infty {\rm e}^{\,-p\, s} 
d_s L^0_s\right),
$$
which in terms of $c_0(\al)$ reads 
$$
\frac 1\pi\int_0^\infty\frac{d\xi}{p+\xi^\al}
=
c_0(\al)\, \int_0^\infty {\rm e}^{\,-p\, s} 
s^{-1/\al}\,ds.
$$
An elementary computation reveals that 
$$
c_0(\al)=\frac 1\pi\,\Gamma((\al+1)/\al),
$$  
hence,
$$
c_1(\al)=\left((\al-1)\,\pi\, m_{\al-1}\right)/\Gamma(1/\al).
$$

Next we find from formula (\ref{e222}) that 
\begin{equation}
\label{e561}
c_6(\al)=1/c_1(\al).
\end{equation}

To compute $c_8(\al,\ga)$ for $\al-1\leq \ga\leq \al/2$ and the limiting case
$c_2(\al)=c_8(\al,\al-1)$ notice from (\ref{e260}) that
$$
c_8(\al,\ga)=\int_\R\nu(dy)\left(|1+y|^\ga -1\right)^2.
$$
Comparing the integral representations of $c_3$ and $c_8$ it is seen
that
\begin{equation}
\label{e56}
2c_3(\al,\ga)+c_8(\al,\ga)=c_3(\al,2\ga)
\end{equation}
which can also be deduced from the following formulae
\begin{eqnarray*}
&&\hskip-1cm(i)\hskip1cm
\E\left(|X_t|^{2\,\ga}\right)
=2\,\E\left(\int_0^t|X_{s}|^{\ga}\,d_sA^{(\ga)}_s\right)
+\E\left(<\hskip-1mm N^{(\ga)}\hskip-1mm>_t\right) \\
&&
\hskip2.5cm
=\left(2c_3(\al,\ga)+c_8(\al,\ga)\right)
\,\E\left(\int_0^t \frac{ds}{|X_s|^{\al-2\ga}}\right),
\\
&&
\hskip-1cm
(ii)\hskip.9cm
\E\left(|X_t|^{2\,\ga}\right)
=c_3(\al,2\ga)\,\E\left(\int_0^t \frac{ds}{|X_s|^{\al-2\ga}}\right).
\end{eqnarray*}   
The first one of these is an easy application of the It\^o formula for
semimartingales and the second one follows (\ref{e52}) because
$\ga\leq \al/2.$ From equation (\ref{e56}) we get 
\begin{eqnarray*}
&&
c_8(\al,\ga)=c_3(\al,2\ga)-2c_3(\al,\ga)
\\
&& 
\hskip1.5cm
=
\frac{2\ga}{\al}\left(\frac{ m_{2\ga}}{m_{2\ga-\al}}-\frac{ m_\ga}{m_{\ga-\al}}\right).\end{eqnarray*} 
The constant $c_2$ is now obtained by letting here $\ga\to\al-1$ 
and using $m_{-1}=+\infty.$ Consequently 
$$
c_2(\al)=\frac{2(\al-1)}{\al}\,\frac{ m_{2(\al-1)}}{m_{\al-2}}
.
$$

To find the constant $c_5(\al),$ we use the relationship (\ref{e201})
between $\Psi$ and $\nu$  which yields
after substitution $y=\xi z$
$$
c_5(\al)=\left(2\,\int_{0}^{\infty}\frac{1-\cos y}{y^{\alpha+1}}\,
dy\right)^{-1}.
$$
Integrating by parts and using
the formulae 2.3.(1) p. 68 in Erdelyi
et al. \cite{erdelyi54} lead us to the explicit value of the integral
$$
\int_{0}^{\infty}\frac{1-\cos y}{y^{\alpha+1}}\,dy
=
\frac{\Gamma(1-\alpha)}{\alpha}\cos(\alpha\pi/2).
$$

The constant $c_6(\alpha)$ can also clearly be expressed in terms of
$c_5$ 
\begin{eqnarray*}
 &&c_6(\alpha)=\left(2\pi c_5(\alpha-1)\right)^{-1}=
\frac 1\pi\,\int_0^\infty\frac{1-\cos \xi}{\xi^{\,\alpha}}\,d\xi\\
&&\hskip1cm
=
\frac 1\pi \frac{\Gamma(2-\alpha)}{\alpha-1}\cos((\alpha-1)\pi/2).
\end{eqnarray*}
It can be verified by the duplication formula for the 
Gamma function  that this agrees with (\ref{e561}).
It holds also that  $c_6(\alpha)\to 1/2$ as $\alpha\uparrow 2.$

The constant $c_7$ is obtained  by simply comparing 
the definitions of $N^x$ in (\ref{c1}) and $\tilde N^x$ in
Proposition \ref{prop1}. We have
$$
N^x_t=\frac{1}{c_6(\al)}\,\tilde N^x_t
$$
implying 
$$
c_7(\alpha)=c_2(\al)\,(c_6(\al))^{\,2}.
$$

\section{Symmetric principal values of local times}
\label{sec6}

Our previous results may be summarized as follows
\begin{enumerate}
\item for $\al-1\leq \ga<\al$ the process $\{|X_t-x|^\ga\}$  is a
  submartingale whose Doob-Meyer decomposition is given by (\ref{e25}),
\item for $(\al-1)/2< \ga<\al-1$ the process $\{|X_t-x|^\ga\}$  is a
  Dirichlet process whose canonical decomposition is given by (\ref{e261}).
\end{enumerate}
\noindent
These results do not discuss whether  $\{(X_t-x)^{\ga,*}\},$ the
symmetric power of order $\ga,$ i.e.,
\begin{equation}
\label{e60}
(X_t-x)^{\ga,*}:={\rm sgn}(X_t-x)\,|X_t-x|^\ga,
\end{equation}
is or is not a semimartingale or a Dirichlet process. In the present
section it is seen that this question can be answered completely
relying on some results in
Fitzsimmons and Getoor \cite{fitzsimmonsgetoor92a} and \cite{fitzsimmonsgetoor92b}, see also
K. Yamada \cite{yamada02}. Let  $x=0$ in (\ref{e60}) and introduce the
principal value integral (cf. (\ref{c4})) 
$$
{\rm p.v.}\,\int_0^t\frac{ds}{X_s^{\theta,*}}
:=\int_0^\infty \frac{dz}{z^\theta }
\left(L^z_t-L^{-z}_t\right),
$$
where by the H\"older continuity (\ref{f30}) the integral is well
defined for $\theta<(\al+1)/2.$ 

\begin{proposition}
\label{prop61}
{\bf a)}
For $\al-1<\ga<\al$ the
process $\{X_t^{\ga,*}\}$ is a semimartingale.\hfill\break\hfill
{\bf b)} For $(\al-1)/2<\ga\leq \al-1$ the
process $\{X_t^{\ga,*}\}$ is a Dirichlet process and not a
semimartingale.
\hfill\break\hfill
{\bf c)} In both cases the unique canonical decomposition of the process
can be written as
\begin{equation}
\label{e61}
X_t^{\ga,*}\, q_{\al,\theta}=N_t^{\ga,*}+c_1(\al)\
{\rm p.v.}\,\int_0^t\frac{ds}{X_s^{\al-\ga,*}},
\end{equation}
where
$$
q_{\al,\al-\ga}=\int_0^\infty\frac
{dx}{x^{\al-\ga }}\left(|1-x|^{\al-1}-(1+x)^{\al-1}\right)
$$
and
$$
N_t^{\ga,*}=\int_0^\infty\frac
{dx}{x^{\al-\ga}}\left(N_t^x-N_t^{-x}
\right).
$$
In particular, for $\ga=\al-1$
\begin{equation}
\label{e62}
X_t^{\al-1,*}\, q_{\al,1}=N_t^{\al-1,*}+c_1(\al)\
{\rm p.v.}\,\int_0^t\frac{ds}{X_s}.
\end{equation}
\end{proposition}

\begin{proof}
Because $\{X_t\}$ is a martingale, it follows from the Ito
formula for semimartingales (\ref{meyer1}) that for $1\leq \ga \leq \al$ the
process $\{X_t^{\ga,*}\}$ is a semimartingale. The other statements in
a) and b) are derived from the decomposition (\ref{e61}) which we now
verify similarly as (\ref{e261}) in Proposition \ref{prop3}. Hence,
we start again from the identity (\ref{c1}) considered at $x$ and
$-x,$ and write, informally 
\begin{eqnarray}
\label{e63}
&&
\int_0^\infty\frac
{dx}{x^{\al-\ga}}\left(|X_t-x|^{\al-1}-|X_t+x|^{\al-1}\right)\\
&&
\nonumber
\hskip1cm
=
\int_0^\infty\frac
{dx}{x^{\al-\ga}}\left(N_t^x-N_t^{-x}
\right)
+c_1(\al)\,\int_0^\infty\frac
{dx}{x^{\al-\ga}}\left(L_t^x-L_t^{-x}\right),
\end{eqnarray}
To analyze the integral on the left hand side consider
$$
Q_{\al,{\al-\ga}}(a)=\int_0^\infty\frac
{dx}{x^{\al-\ga}}\left(|a-x|^{\al-1}-|a+x|^{\al-1}\right).
$$
It is easily seen that this integral is absolutely convergent and 
$$
Q_{\al,{\al-\ga}}(a)=a^{\ga,*}\,q_{\al,{\al-\ga}}.
$$
Now the rest of the proof  is very similar to that of Proposition
\ref{prop3} b), and  is therefore omitted. 
\end{proof}
\begin{remark}
{\bf a)}
{\rm  The increasing process associated with $N^{\ga,*}$ is given by
\begin{eqnarray*}
&&
<\hskip-1mm N_t^{\ga,*}\hskip-1mm >_t\, =
q_{\al,\al-\ga}^2\,
 \int_0^tds\int_\R \nu(dz)\left((X_s+z)^{\al-\ga,*}-X_s^{\al-\ga,*}\right)^2\\
 &&
 \hskip1.6cm
 =
q_{\al,\al-\ga}^2\,
 \int_0^t\frac{ds}{|X_s|^{\al-2\ga}}\,\int_\R \nu(dz)\left((1+z)^{\al-\ga,*}-1\right)^2.
\end{eqnarray*}
We also have by scaling 
\begin{eqnarray*}
&&
\E\left(\int_0^t\frac{ds}{|X_s|^{\al-2\ga}}\right)
=\int_0^t\,s^{(2\ga-\al)/\al}\,ds\, \E\left(
|X_1|^{2\ga-\al}\right)\\ 
&&
\hskip3.3cm
=\frac{\al}{2\ga}\ t^{2\ga/\al}
\, \E\left(
|X_1|^{2\ga-\al}\right).
\end{eqnarray*}
}
\hfill\break\hfill
{\bf b)} {\rm Since
$$
|X_t|^\gamma= (X_t^+)^\gamma + (X_t^-)^\gamma 
$$
and 
$$
|X_t|^{\gamma,*}= (X_t^+)^\gamma - (X_t^-)^\gamma 
$$
it is straightforward to derive the decomposition formulae for
$\{(X_t^+)^\gamma\}$ and $\{ (X_t^-)^\gamma\}$, and we leave this to
the reader.
}
\hfill\break\hfill
{\bf c)} {\rm  Note how different (\ref{e62}) is in the Brownian case
   $\al=2,$ for which on one hand $\{B_t\}$ is a martingale, and on
   the other hand
$$
\varphi(B_t)=\int_0^t\log|B_s|\,dB_s + \frac 12\ {\rm p.v.}\,\int_0^t\frac
    {ds}{B_s}
$$
with $\varphi(x)=x\log|x|-x.$ For principal values of Brownian motion
and extensions of It\^o's formula, see Yor \cite{yor82}, \cite{yor97} and  Cherny \cite{cherny01}.
}
\end{remark}

\bibliographystyle{plain}
\bibliography{yor1}
\end{document}